\documentclass[12pt,a4paper]{article}%
\usepackage{graphicx}
\usepackage{amsmath}%
\usepackage{amsfonts}%
\usepackage{amssymb}
\newtheorem{theorem}{Theorem}

\newtheorem{thm}{Theorem}[section]
\newtheorem{cor}{Corollary}[section]
\newtheorem{cla}[cor]{Claim}
\newtheorem{lem}[cor]{Lemma}
\newtheorem{prop}[cor]{Proposition}
\newtheorem{defe}{Definition}
\newenvironment{proof}[1][Proof]{\textbf{#1.} }{\ \rule{0.5em}{0.5em}}
\begin{document}

\title{When is $ch(K_{m,n}) = m-1$?}
\author{Nurit Gazit} \date{June 2006}
\maketitle
\thispagestyle{empty} %

\newpage

\begin{abstract}
Let $n_m$ be the smallest integer $n$ such that $ch(K_{m,n}) = m-1$,
where $ch(G)$ denotes the choice (list chromatic) number of the
graph $G$. We prove that there is an infinite sequence of integers
$S$, such that if $m \in S$, then $n_m \leq 0.4643 {(m-2)}^{m-2}$.
If $m \rightarrow \infty$, then $n_m$ is asymptotically at most
$0.474 {(m-2)}^{m-2}$.

\end{abstract}

\newpage

\section{Introduction}


A \emph{list assignment} of a graph $G = (V,E)$ from a family of
sets (color lists) $L$ is an assignment to each vertex $v \in V(G)$
of a list $L(v)$ of colors. A \emph{k-list assignment} is a list
assignment that satisfies $|L(v)| = k$ for every $v \in V(G)$. An
\emph{$L$-coloring} is a function $c: V(G) \rightarrow \bigcup_{v
\in V(G)}{L(v)}$ that assigns each vertex $v$ a color $c(v) \in
L(v)$. A \emph{proper $L$-coloring} is an $L$-coloring such that the
neighbors of each vertex $v$ are colored in a different color than
that of $v$. A graph $G$ for which there is a proper $L$-list
coloring is called \emph{$L$-list choosable}. The \emph{choice
number} $ch(G)$ is the minimum number $k$ such that for every k-list
assignment $L$ of $G$, there is a proper $L$-coloring of $G$. The
concept of choosability was introduced by Vizing in 1976 \cite{[6]}
and independently by Erd\H{o}s, Rubin and Taylor in 1979 \cite{[1]}.
It is also shown in \cite{[1]} that the choice number of the
complete bipartite graph $K_{n,n}$ satisfies $ch(K_{n,n}) = (1 +
o(1)) \log_2{n}$.

In \cite{[2]}, Gazit and Krivelevich calculate the asymptotic value
of the choice number of complete multi-partite graphs where the
sizes of the different parts are not too far apart, i.e. of graphs
of the form $K_{n_0,...,n_s}$, $n_0 \leq n_1 ... \leq n_s$, where
$n_0$ is not too small compared to $n_s$.

In particular, for the bi-partite case, Gazit and Krivelevich prove:
\begin{itemize}
\item[($\bullet$)]
Let $2 \leq n_0 \leq n_1$ be integers, and let $n_0 =
{(\log{n_1})}^{\omega(1)}$. Denote ${ k = \frac{ \log{n_1}}{
\log{n_0}} }$. Let $x_0$ be the unique root of the equation $x - 1 -
x^{\frac{k-1}{k}} = 0$ in the interval $[1, \infty)$. Then
$ch(K_{n_0,n_1}) = (1+o(1))\frac{\log{n_1}}{ \log{x_0}}$.
\end{itemize}

While the above-mentioned result deals with bipartite graphs
$K_{m,n}$ in which $m$ is not too small compared to $n$, in this
paper we consider bipartite graphs $K_{m,n}$ in which $m$ is very
small compared to $n$.

It is trivial to see that $ch(K_{m,n}) = m+1$ if $n \geq m^m$. To
see this, let $(M,N)$ $(|M| = m, |N| = m^m)$ be a bipartition of
$K_{m,m^m}$. Assign $m$ pairwise-disjoint lists of colors to $M$,
and assign all $m^m$ colorings of these $m$ lists as color lists to
the vertices of $N$. Clearly, $K_{m,m^m}$ is not choosable from
these lists. This shows that $ch{(K_{m,m^m})} > m$ (and therefore
also $ch{(K_{m,n})} > m$ for $n \geq m^m$). But $ch(K_{m,n}) \leq
m+1$ for every $n$, since for every $(m+1)$-list assignment $L$,
every coloring of the vertices on $M$ uses at most $m$ different
colors, leaving at least one color $d(v) \in L(v)$ in the color list
of every $v \in N$ which has not been used - assigning $d(v)$ to
every $v \in N$, completes a proper $L$-coloring.

Hoffman and Johnson \cite{[3]} proved that $ch(K_{m,n}) = m$ \\
if and only if ${{(m-1)}^{m-1} - {(m-2)}^{m-1}} \leq n < m^m$.

Let $n_m$ be the smallest integer $n$ such that $ch(K_{m,n}) = m-1$.
In this paper we aim to find an upper bound on $n_m$. We will show:\\
\textbf{Theorem 1} \textsl{If $n_m = p {(m-2)}^{m-2}$, and
$m'-2=k(m-2)$, with $k$ integer, then $n_{m'} \leq p
{(m'-2)}^{m'-2}$.}
\newline
\textbf{Theorem 2} \textsl{There is an infinite sequence of integers
$S$, such that if $m \in S$, then $n_m \leq 0.4643 {(m-2)}^{m-2}$.}
\newline
\textbf{Theorem 3} \textsl{If $m \rightarrow \infty$, then $n_m$ is
asymptotically at most $0.474 {(m-2)}^{m-2}$.}

\newpage

\section{Definitions and Preliminary Observations}


We begin by giving several definitions. Note that in this paper all
sets are of finite cardinality.

A \emph{hypergraph} $H$ is a pair $H = (V, E)$, where $V$ is a set
of elements called \emph{vertices}, and $E$ is a set of subsets of
$V$ called \emph{hyperedges}.

A \emph{transversal} of a family of sets $S$ is a set $S_t$ such
that $s \cap S_t \neq \emptyset$ for all $s \in S$.

The \emph{transversal hypergraph} of a family of sets $S$, which
will be denoted by $T(S)$ is the hypergraph whose vertices are the
union of the sets in $S$, and whose edges are all the transversals
of cardinality $\leq |S|$ of $S$. The transversal hypergraph of the
set of edges of a hypergraph $H$ will be denoted by $T(H)$.


Given an ordered family of sets $L$ of cardinality $m$, a
\emph{track} is an ordered $m$-tuple, created by choosing element
$c_i$ from set $L_i$, $1 \leq i \leq m$.

Given a family of sets $S$ and a transversal $e$ of $S$, a track $t$
\emph{belongs} to $e$ if the set of distinct elements of $t$ is $e$.
We say that the transversal $e$ \emph{represents} the track $t$.


Let $R = (V_R, E_R)$, $H = (V_H, E_H)$ be hypergraphs. An
\emph{$R$-cover} of $H$ is a sub-hypergraph $C$ of $R$ such that
every hyperedge of $H$ contains at least one hyperedge $c \in E_C$.

A \emph{$k$-cover} of a hypergraph $H$ is an $R$-cover of $H$, $R$
being all the subsets of size $k$ of the vertices of $H$ (though the
subsets can also be taken from a larger set).

A \emph{minimum $R$-cover} is an $R$-cover whose edge set has the
least cardinality among those of all $R$-covers. The cardinality of
a minimum $R$-cover of $H$ will be called \emph{$cov_{R}(H)$}. If
$H$ allows no $R$-cover, $cov_{R}(H) = \infty$.

A \emph{minimal $R$-cover} is an $R$-cover which does not contain
any other $R$-cover.

Clearly, a minimum $R$-cover is also a minimal $R$-cover. It is easy
to see that the following lemma holds.

\begin{lem}\label{le1}
Let $R'$ be a minimal $R$-cover of a hypergraph $H$ with $E_H \neq
\emptyset$. If $e_R \in R'$, then there is at least one edge $e_H
\in E_H$ such that $e_R$ is a unique edge of $R'$ which is a subset
of $e_H$.

\end{lem}



A \emph{minimal edge} of a hypergraph $H = (V,E)$ is an edge that
does not contain any other edge in $H$.

Since the edges are of finite cardinality, every edge $e$ of $H$
contains a minimal edge (for example, an edge $h$ of $H$ of minimum
cardinality contained in $e$, is obviously a minimal edge).

Removing edges from a hypergraph cannot turn a minimal edge into a
non-minimal edge. Removing non-minimal edges from a hypergraph
cannot turn a non-minimal edge to a minimal edge, since every
non-minimal edge contains a minimal edge.

Therefore removing non-minimal edges from a hypergraph does not
change the set of minimal edges of the hypergraph.


A \emph{minimal hypergraph} is a hypergraph in which every edge is
minimal (i.e. a graph in which no edge is contained in any other
edge).

The \emph{cover hypergraph} of a hypergraph $H$ is the
sub-hypergraph of $H$ whose edge set is all the minimal edges of
$H$.

The cover hypergraph is a minimal hypergraph, since removing edges
from a hypergraph can not change a minimal edge to a non-minimal
one.

Vertices $a$ and $b$, $a \neq b$ of a hypergraph $H$ are called
\emph{min-equal} if:

\begin{defe} for every minimal edge $e \in E_H$, if $a \in e$ then
there is a minimal edge $e' \in E_H$ such that $a \notin e'$, $b \in
e'$, and $v \in e' \Rightarrow v \in e$, and likewise if $b \in e$.
\end{defe}

\begin{defe} for every minimal edge $e \in E_H$, if $a \in e$ then $e' = e
\setminus \{ a \} \cup \{ b \}$ is a minimal edge of $H$, and
likewise if $b \in e$.
\end{defe}

It is easy to see the two definitions are equivalent.

A vertex $a$ is min-equal to itself.

It is easy to see that min-equality is an equivalence relation.

\begin{lem}\label{lemin}
If $a$ and $b$ are min-equal vertices of $H$, $a \neq b$, then no
minimal edge contains both $a$ and $b$.
\end{lem}
\begin{proof} If $e \in H$ is minimal and $a \in e$, there is a
minimal edge $e' \in E_H$ such that $a \notin e'$, $b \in e'$, and
$v \in e' \Rightarrow v \in e$.  But then if $b \in e$, $e'
\subsetneq e$ ($a \in e$ but $a \notin e'$) in contradiction to the
minimality of $e$.
\end{proof}
\bigskip

Vertices $a$ and $b$ of a hypergraph $H$ are called \emph{H-equal},
if for every edge $e \in E_H$, if $a \in e$ and $b \notin e$, then
$e' = e \setminus \{ a \} \cup \{ b \}$ is also a hyperedge of $H$,
and likewise if $b \in e$.

\section{Proof of Theorem 1}




\begin{theorem}\label{thm1}
Let integers $m$, $m'$ satisfy
$m'-2 = k(m-2)$ for some integer k. If $n_m = p {(m-2)}^{m-2}$ for
some real $p$, then $n_{m'} \leq p {(m'-2)}^{m'-2}$.
\end{theorem}

To prove this theorem, we shall first restate the problem of finding
$n_m$ in terms of finding minimum $R$-covers of hypergraphs.

Hoffman and Johnson \cite{[3]} prove the following simple lemma:

\begin{lem}\label{le0} If $L$ is a list assignment to $K_{m,n}$, with
bipartition $(M,N)$ $(|M| = m, |N| = n)$, then there is no proper
$L$-coloring of $K_{m,n}$ if and only if each transversal of the
sets $L(v)$, $v \in M$, contains one of the sets $L(u)$, $u \in N$.
\end{lem}

\begin{cor}
Let $R_m$ be the hypergraph whose vertices are \\ $S =
\{1,2,...,m(m-2)\}$ and whose edges are all subsets of cardinality
$m-2$ of S. Then $n_m$ is equal to the minimum over all
sub-hypergraphs $R'$ of $R_m$ with $m$ edges, of $cov_{R_m}{T(R')}$.
\end{cor}

\begin{proof}
Let $min\_{cov}(m)$ be the minimum over all sub-hypergraphs $R'$ of
$R_m$ with $m$ edges, of $cov_{R_m}{T(R')}$. Let us take a
sub-hypergraph with $m$ edges $R'$ and a minimum $R_m$-cover $R''$
of its transversal hypergraph $T(R')$, such that $|R''| =
min\_{cov}(m)$. Then assigning the edges of $R'$ (a family of $m$
sets of size $m-2$) to the vertices of $M$, and the edges of $R''$
(a family of $min\_{cov}(m)$ sets of cardinality $m-2$) to the
vertices of $N$, gives, by Lemma \ref{le0}, an $(m-2)$-list
assignment $L$ of $K_{m,min\_{cov}(m)}$ such that there is no proper
$L$-coloring of $K_{m,min\_{cov}(m)}$, implying that
$K_{m,min\_{cov}(m)}$ is not $(m-2)$-choosable, and $n_m \leq
min\_{cov}(m)$.

By the definition of $n_m$, there is an $(m-2)$-list assignment $L$
to $K_{m,n_m}$ such that there is no proper $L$-coloring of
$K_{m,n_m}$. The color lists assigned to $M$ in $L$ are $m$ lists of
cardinality $m-2$, so there are at most $m(m-2)$ colors in their
union $\bigcup_{v \in M}{L(v)}$, and we may label them WLOG as
$\{1,2,...m(m-2)\}$ (therefore they make up a sub-hypergraph $R'$ of
$R_m$ with $m$ edges).  If any of the lists $L(u)$, $u \in N$
contains a color other than $\{1,2,...,m(m-2)\}$, then no
transversal of the sets $L(v)$, $v \in M$ contains $L(v)$. Therefore
if we remove $v$ (and $L(v)$) from $N$ to get $N'$, $|N'| =
{n_m}-1$,
Lemma \ref{le0} gives us $ch(K_{m,n_m - 1}) > m-2$, in contradiction
to the minimality of $n_m$. Therefore the color lists assigned to
$N$ are $n_m$ subsets of cardinality $m-2$ of $\{1,2,...,m(m-2)\}$ -
i.e. a sub-hypergraph $R''$ of $R_m$ of cardinality $n_m$. By Lemma
\ref{le0} each transversal of the sets $L(v)$, $v \in M$, contains
one of the sets $L(u)$, $u \in N$, i.e. $R''$ is an $R_m$-cover of
$T(R')$, therefore $n_m \geq cov_{R_m}{T(R')} \geq min\_{cov}(m)$.

\end{proof}
\bigskip

Now we shall show that when looking for a minimum $R$-cover of a
hypergraph, it suffices to look for a minimum $R$-cover of its cover
hypergraph.

It is easy to see that the following lemma holds:

\begin{lem}
Let $H$, $H'$ and $R$ be hypergraphs. If every $R$-cover of $H'$ is
an $R$-cover of $H$, and if $R'$ is a minimum (minimal) $R$-cover of
$H$ and it is also an $R$-cover of $H'$, then $R'$ is a minimum
(minimal) $R$-cover of $H'$.
\end{lem}


\begin{cor}\label{co0}
Let $H$, $H'$ and $R$ be hypergraphs. If for every sub-hypergraph
$R'$ of $R$, $R'$ is an $R$-cover of $H$ if and only if $R'$ is an
$R$-cover of $H'$, then $R'$ is a minimum (minimal) $R$-cover of $H$
if and only if $R'$ is a minimum (minimal) $R$-cover of $H'$.
\end{cor}



\begin{lem}\label{le2}
Given a hypergraph $H$, let $H_{cov}$ be the cover hypergraph of
$H$. Let $R$ be another hypergraph. Then $R'$ is an $R$-cover of $H$
if and only if $R'$ is an $R$-cover of $H_{cov}$. Therefore, by
Corollary \ref{co0}, $R'$ is a minimum (minimal) $R$-cover of $H$ if
and only if $R'$ is a minimum (minimal) $R$-cover of $H_{cov}$.
\end{lem}

\begin{proof}
If $R'$ is an $R$-cover of $H$ it is also an $R$-cover of $H_{cov}$,
since $H_{cov}$ is a sub-hypergraph of $H$. Now let $R'$ be an
$R$-cover of $H_{cov}$. We will show that for every $e \in H$ there
is an $r \in R'$ which $e$ contains. If $e$ is a minimal edge then
$e \in H_{cov}$, and since $R'$ is an $R$-cover of $H_{cov}$, there
is an $r \in R'$ which $e$ contains. If $e$ is a non-minimal edge
then $e$ contains a minimal edge $e'$. Since $e'$ is minimal it is
in $H_{cov}$, therefore there is an $r' \in R'$ so that $r'
\subseteq e' \subseteq e$.
\end{proof}

\begin{prop}\label{the1}
Let $R$ and $H$ be hypergraphs. Let $\{ V_\alpha \}, \alpha \in A$,
be a partition of $H$ into disjoint sets of min-equal, $R$-equal
vertices.
For every $\alpha \in A$ let us choose a representative $v_\alpha
\in V_\alpha$, and attach a weight $w(v_\alpha) = |V_\alpha|$ to it.
Let $H_{cov}(\{ v_\alpha \}_{\alpha \in A})$ be the sub-hypergraph
of $H$ whose edges are the minimal edges of $H$ contained in
$\{v_\alpha\}_{\alpha \in A}$. Then $cov_{R}(H)$ is equal to the
minimum over all $R$-covers $R'$ of $H_{cov}(\{ v_\alpha \}_{\alpha
\in A})$ of $\sum_{r \in R'}{\{\prod_{v \in r} w(v) \} }$.
\end{prop}

\begin{proof}
By Lemma \ref{le2}, $R'$ is a minimum $R$-cover of $H$ if and only
if $R'$ is a minimum $R$-cover of $H_{cov}$.
\begin{cla}
Let $R'$ be a minimum $R$-cover of $H_{cov}$, and let $\alpha \in
A$. Then $R'$ partitions the edges of $H_{cov}$ into the following
$|V_\alpha| + 1$ sets:\\ $W_{\emptyset}$ - edges that contain an
$e_R \in R'$ such
that $V_\alpha \cap e_R = \emptyset$; \\
For each $v_{\beta} \in V_\alpha$, a set $W_{\beta}$ such that if
$e_H \in W_{\beta}$, then $v_{\beta} \in e_H$, and $V_\alpha \cap
e_R = \{v_{\beta}\}$ for
every $e_R \in R'$ s.t. $e_R \subseteq e_H$.\\
In this partition, for every $v_{\beta}, v_{\beta'} \in V_{\alpha}$,
$e_H \in W_{\beta} \Rightarrow e_H \setminus \{ v_{\beta} \} \cup \{
v_{\beta'} \} \in W_{\beta'}$. 
\end{cla}
\begin{proof}
Let $e_H$ be an edge in $H_{cov}$. Then since $e_H$ is minimal and
the vertices in $V_\alpha$ are min-equal, by Lemma \ref{lemin}
either $V_\alpha \cap e_H = \emptyset$ or there is a $v_\beta \in
V_\alpha$ such that $V_\alpha \cap e_H = \{v_\beta\}$. If $V_\alpha
\cap e_H =
\emptyset$ then $e_H \in W_{\emptyset}$ ($R'$ is a cover). \\
Suppose $V_\alpha \cap e_H = \{v_\beta\}$. If there is an $e_R \in
R'$, $e_R \subseteq e_H$ such that $V_\alpha \cap e_R = \emptyset$,
then $e_H \in W_{\emptyset}$. Otherwise, $e_H \in W_{\beta}$. This
gives us the partition. $e_H \in W_{\beta} \Rightarrow e_H \setminus
\{ v_{\beta} \} \cup \{ v_{\beta'} \} \in W_{\beta'}$ because since
$v_{\beta}$ and $v_{\beta'}$ are min-equal, if $e_H \in H_{cov}$,
then $e_H' = e_H \setminus \{ v_{\beta} \} \cup \{ v_{\beta'} \} \in
H_{cov}$. Therefore either $e_H' \in W_{\emptyset}$ or $e_H' \in
W_{{\beta}'}$. But if there is an $e_R \in R'$, $e_R \subseteq e_H'$
such that $V_\alpha \cap e_R = \emptyset$, then since $v_\beta'
\notin e_R$ and $e_R \subseteq e_H'$, $e_R \subseteq e_H$, in
contradiction to the assumption that $e_H$ contains no such $e_R$.
\end{proof}

\begin{cla}\label{cl2}
Let $v_{\beta} \in V_{\alpha}$. Let $R_{\beta}$ be the edges of $R$
that contain $v_{\beta}$ and do not contain any other $v_{\beta'}
\in V_\alpha$. Then $R'$ induces a minimum $R_{\beta}$-cover
${R_{\beta}}*$ on $W_{\beta}$, i.e. the set of $e_R \in R'$ s.t.
$e_R \subseteq e_H$ for some $e_H \in W_{\beta}$ is a minimum
$R_{\beta}$-cover of $W_{\beta}$.
\end{cla}
\begin{proof}
The edges in $R'$ used to cover $W_{\beta}$ (i.e. those contained in
edges of $W_{\beta}$) all contain $v_{\beta}$ and not any other
$v_{\beta'} \in V_\alpha$, so the induced cover is an
$R_{\beta}$-cover. If it is not a minimum $R_{ \beta}$-cover, then
we can replace the edges used to cover $W_{\beta}$ with those of a
minimum $R_{\beta}$-cover, and get an $R$-cover of smaller
cardinality of $H_{cov}$ (every edge not in $W_{ \beta}$ contains an
edge in $R'$ that does not contain $v_{\beta}$ - so this is still an
$R$-cover), in contradiction to the fact that $R'$ is a minimum
$R$-cover.
\end{proof}
\bigskip

Replacing the $R_{\beta}$-cover on $W_{ \beta}$ with another minimum
$R_{\beta}$-cover without touching the $r \in R'$ which do not
contain $v_{\beta}$ still gives a minimum $R$-cover (since it is
still a cover - we use only edges in R, and it is of the same
cardinality as the previous $R$-cover).


If $v_{\beta}$ and $v_{\beta'}$ are $R$-equal, then the edges in
$R_{\beta}$ are the edges in $R_{\beta}'$ with $v_{\beta}$ changed
to $v_{\beta'}$.  Therefore, since there is an isomorphism between
$W_{ \beta}$ and $W_{{\beta'}}$ and between $R_{ \beta}$ and $R_{
{\beta'}}$ in which $v_{\beta} \rightarrow v_{{\beta'}}$, and the
other elements stay the same, the minimum $R_{ \beta}$-covers are
also isomorphic, and therefore are of the same cardinality. Thus, if
we take the $R_{\beta}$-cover of $W_{\beta}$ used in $R'$, then for
every $v_{{\beta'}} \in V_{\alpha}$, replacing $v_{\beta}$ in every
edge of the $R_{\beta}$-cover by $v_{{\beta'}}$ gives a minimum
$R_{{\beta'}}$-cover of $W_{{\beta'}}$, and thus we get a new
minimum $R$-cover $R''$.

If $e_H \in W_\emptyset$ and there is an $e_R \in R''$ such that
$e_R \subseteq e_H$ and $e_R \cap V_\alpha = \{v_{\beta}\}$, then
there is an $e_H' \in W_\beta$ such that $e_R \subseteq e_H'$
(otherwise we can remove $e_R$ from $R''$ and still get a cover).
Therefore, for every $v_{\beta} \in V_{\alpha}$, the cardinality of
$R''$ is equal to $| \{ e_R \in R'' | e_R \cap V_{\alpha} =
\emptyset \} |$ + $|V_{\alpha}| |\{ e_R \in R'' | e_R \cap
V_{\alpha} = \{v_{\beta}\} \} |$. This is because according to the
previous statement, for every $v_{{\beta}'} \in V_{\alpha}$, if
$e_R' \in R''$, then $e_R' \cap V_{\alpha} = \{v_{\beta'}\}$ if and
only if there is an $e_H \in W_{{\beta}'}$ such that $e_R' \subseteq
e_H$, i.e. if and only if $e_R'$ belongs to the $R_{{\beta}'}$-cover
of $W_{{\beta}'}$ induced by $R''$. So $| \{ e_R \in R'' | e_R \cap
V_{\alpha} = \{v_{{\beta}'}\} \} |$ is equal to the cardinality of
the $R_{{\beta}'}$-cover of $W_{{\beta}'}$ induced by $R''$, and
these are all equal to the cardinality of the $R_{{\beta}}$-cover of
$W_{{\beta}}$ induced by $R''$, which is equal to $| \{ e_R \in R''
| e_R \cap V_{\alpha} = \{v_{{\beta}}\} \} |$.

Let us define $R_\emptyset$ to be the sub-hypergraph of $R$ whose
edges satisfy ${e_R \cap V_{\alpha} = \emptyset}$. Given a
$v_{\beta} \in V_\alpha$, let us define $H_{cov}{\{v_{\beta}\}}$ to
be the sub-hypergraph of $H_{cov}$ composed of all edges $e_H \in
H_{cov}$ such that $e_H \cap V_{\alpha} \subseteq \{v_{\beta}\}$.

Let $R'' \{ v_{\beta} \}$ be the $\{ R_{\emptyset} \cup R_{\beta}
\}$-cover that $R''$ induces on $H_{cov}{\{v_{\beta}\}}$.

Observe that $\{e_R \in R'' \{ v_{\beta} \} | e_R \cap V_{\alpha} =
\emptyset \} = \{ e_R \in R'' | e_R \cap V_{\alpha} = \emptyset \}$.
This is because if $e_R \in R''$ and $e_R \cap V_{\alpha} =
\emptyset$, then since $R''$ is minimal, there is an $e_H \in
W_{\emptyset}$ such that $e_R \subseteq e_H$. If $e_H \in
H_{cov}{\{v_{\beta}\}}$ then $e_R \in R'' \{ v_{\beta} \}$.
Otherwise, $e_H \cap V_{\alpha} = \{v_{\beta'}\}$, and since $e_R
\cap V_{\alpha} = \emptyset$, $e_R \subseteq e_H \setminus \{
v_{\beta'} \} \cup \{ v_{\beta} \} \in H_{cov}{\{v_{\beta}\}}$, so
$e_R \in R'' \{ v_{\beta} \}$.

Also $\{ e_R \in R'' \{ v_{\beta} \} | e_R \cap V_{\alpha} =
\{v_{\beta}\} \} = \{ e_R \in R'' | e_R \cap V_{\alpha} =
\{v_{\beta}\} \}$. This is because if $e_R \in R''$ and $e_R \cap
V_{\alpha} =\{ v_{\beta} \}$ then since $R''$ is minimal there is an
$e_H$ such that $e_R \subseteq e_H$, but then $e_H \in
H_{cov}{\{v_{\beta}\}}$ and $e_R \in R_{\beta}$, so $e_R \in R'' \{
v_{\beta} \}$.

Therefore, there is an $\{ R_{\emptyset} \cup R_{\beta} \}$-cover
$R^{*}$ of $H_{cov}{\{v_{\beta}\}}$ such that $cov_R{(H)} = |\{ r
\in R^{*} | v_{\beta} \notin r\}| + | V_\alpha | \times |\{ r \in
R^{*} | v_{\beta} \in r\}|$.

In order to find a minimum $R$-cover of $H_{cov}$ it suffices to
find an $\{ R_{\emptyset} \cup R_{\beta} \}$-cover $R^{**}$ of
$H_{cov}{\{v_{\beta}\}}$, in which $|\{ r \in R^{**} | v_{\beta}
\notin r\}| + | V_\alpha | \times |\{ r \in R^{**} | v_{\beta} \in
r\}|$ is minimal. This is equivalent to putting a weight of
$w(v_{\beta}) = |V_\alpha|$ on $v_{\beta}$ and 1 on all other
vertices in $H$, and finding a minimum weighted cover - i.e. each $r
\in R$ is given a value of $\prod_{i \in r} {w(i)}$, and we want to
find a cover for which the sum of the values of the edges in the
cover is minimal.

To get an $R$-cover of $H_{cov}$ from $R^{**}$, we simply take the
edges in $R^{**}$ that contain $v_{\beta}$ and add edges in which it
is replaced by $v_{\beta'}$, for every ${v_{\beta'} \in V_\alpha}$.
This is an $R$-cover of $H_{cov}$ because if $e_H \in \{H_{cov}
\setminus H_{cov}{\{v_{\beta}\}} \}$, then ${e_H \cap V_{\alpha} =
\{ v_{\beta'} \}}$, so $e_H' = e_H \setminus \{ v_{\beta'} \} \cup
\{ v_{\beta} \} \in H_{cov}{\{v_{\beta}\}}$. If there is an edge
$e_R$ in $R^{**}$, $e_R \subseteq e_H$ such that $e_R \cap
V_{\alpha} = \emptyset$, then $e_R \subseteq e_H'$. Otherwise there
is an edge $e_R$ in $R^{**}$, $e_R \subseteq e_H$ such that $e_R
\cap V_{\alpha} = v_{\beta}$, so $e_R \setminus \{ v_{\beta} \} \cup
\{ v_{\beta'} \} \subseteq e_H'$.

The cardinality of this $R$-cover will be $|\{ r \in R^{**} |
v_{\beta} \notin r\}| + {| V_\alpha | \times |\{ r \in R^{**} |
v_{\beta} \in r\}| \leq |\{ r \in R^{*} | v_{\beta} \notin r\}|} + |
V_\alpha | \times |\{ r \in R^{*} | v_{\beta} \in r\}| =
cov_R{(H)}$, and therefore it is a minimum $R$-cover.

Now given $v_{\beta} \in V_{\alpha}$ we have a problem of finding a
minimum weighted $\{ R_{\emptyset} \cup R_{\beta} \}$-cover $R'$ of
$H_{cov}{\{v_{\beta}\}}$. Given $V_{\alpha'}$, we can repeat this
process. Given a minimum $\{ R_{\emptyset} \cup R_{\beta} \}$-cover
$R'$, we divide the edges of $H_{cov}{\{v_{\beta}\}}$ into those
that contain $r \in R'$ that are disjoint from $V_{\alpha'}$ - some
of which contain only $r \in R'$ that intersect $V_\alpha$, and
those that do not contain any such $r \in R'$, which are divided
into $|V_{\alpha'}|$ isomorphic sets $W_\beta'$ that contain only
one $v_{\beta'} \in V_{\alpha'}$ (and the isomorphism changes only
the $v_\beta$'s). The same proof shows that this problem is
identical to that in which for one $v_{\beta'} \in V_{\alpha'}$ we
look only at edges $e \in H_{cov}{\{v_{\beta}\}}$ such that $e \cap
V_{\alpha'} \subseteq \{ v_{\beta'} \}$, put a weight on
$v_{\beta'}$ of $|V_{\alpha'}|$, and look for a minimum weighted
cover.

If we continue this process with all the $V_\alpha$'s, we get that
to find a minimum $R$-cover of $H$, it is enough to choose one
representative $v_\alpha$ for each $V_\alpha$, and look only at the
sub-hypergraph $H' = H_{cov}(\{ v_\alpha \}_{\alpha \in A})$ of $H$
which is composed of all the minimal edges of $H$ contained in
$\{v_\alpha\}_{\alpha \in A}$. Each representative $v_\alpha$ is
assigned a weight of $w(v_\alpha) = |V_\alpha|$, and we look for an
$R$-cover $R'$ of $H'$ which minimizes the expression $\sum_{r \in
R'}{\{\prod_{v \in r} w(v) \} }$.

This ends the proof of Proposition \ref{the1}.

\end{proof}

Now we have only two small lemmas left to prove before we reach our
theorem.

\begin{lem}\label{leb4} Let $H$ be the transversal hypergraph $T(L)$
of a family of lists $L$. Let $v_1$ be a vertex in $H$ that appears
in a subset $L'$ of the family of lists $L$. Then every minimal
transversal in $H$ containing $v_1$ represents at least one track in
which $v_1$ is chosen out of every $l \in L'$.
\end{lem}

\begin{proof}
If a transversal $e$ contains $v_1$, then in every track that $e$
represents, $v_1$ is chosen out of some list $l \in L'$. Let us take
such a track $t$, and create a track $t'$ by changing the element
chosen out of every list $l \in L'$ to $v_1$. The transversal that
represents $t'$ is contained in $e$ (we did not add elements to the
transversal), so since the transversal $e$ is minimal, $t'$ belongs
to $e$.
\end{proof}

\begin{lem}\label{le4} Let $H = T(L)$. Let $v_1, v_2$ be vertices in $H$.
Assume that for every list $l \in L$, $v_1 \in l$ if and only if
$v_2 \in l$. Then $v_1$ and $v_2$ are min-equal in $H$.
\end{lem}

\begin{proof}

Let $L'$ be the subset of $L$ of lists that contain $v_1$ and $v_2$.
Let $e \in H$ be a minimal edge in $H$, such that, WLOG, $v_1 \in
e$. Then by Lemma \ref{leb4} there is a track $t$ that belongs to
$e$ in which $v_1$ is chosen out of every $l \in L'$. Let us take
the track $t'$ in which we choose the same color we chose in $t$ out
of every $l \in L \setminus L'$, and choose $v_2$ out of every $l
\in L'$. Then the transversal of the new track is $e' = e \setminus
\{ v_1 \} \cup \{ v_2 \}$.

Suppose $e'$ is not minimal. If $e'' \subseteq e'$, then $v_2 \in
e''$ (otherwise $e'' \subsetneq e$). Take a minimal edge $e''
\subseteq e'$. Since $e'' \neq e'$, $v_2 \in e''$, there is a $v \in
e'$, $v \neq v_2$ (thus $v \in e$), s.t. $v \notin e''$. Since $v_2
\in e''$, by Lemma \ref{leb4} $e''$ represents at least one track
$t''$ in which $v_2$ is chosen out of every $l \in L'$. Creating a
new track which is the same as $t''$ except we choose $v_1$ out of
every $l \in L'$, gives a transversal which is contained in $e$ but
not equal to it, in contradiction to $e$ being minimal.

\end{proof}

We are now ready to prove Theorem \ref{thm1}. First let us
re-formulate it.\\ \textbf{Theorem \ref{thm1}} \textit{Let integers
$m$, $m'$ satisfy $m'-2 = k(m-2)$ for some integer k. Then, for
every family $L_0$ of $m$ ${(m-2)}$-tuples (i.e. a sub-hypergraph of
cardinality $m$ of $R_m$) such that the cardinality of a minimum
$(m-2)$-cover of $T(L_0)$ is $s {(m-2)}^{m-2}$ (for $s$ real), there
is a family $L_0'$ of $m'$ ${(m'-2)}$-tuples such that the
cardinality of a minimum $(m'-2)$-cover of $T(L_0')$ is at most $s
{(m'-2)}^{m'-2}$. Therefore, if $n_m = p {(m-2)}^{m-2}$ for some
real $p$, then $n_{m'} \leq p {(m'-2)}^{m'-2}$.}

\begin{proof}
Given a family $L_0$ of $m$ ${(m-2)}$-tuples $l_i$ $1 \leq i \leq
m$, we create a family $L_0'$ of $m'$ ${(m'-2)}$-tuples $l_i$ $1
\leq i \leq m$, as follows - we take $k$ families of
${(m-2)}$-tuples $L_j$ ($1 \leq j \leq k$) isomorphic to $L_0$, each
in colors which have not been used by the previous families. Let us
denote the $m$ members of $L_j$ by $l_{ji}$, $1 \leq i \leq m$.
We put the lists $L_j$ side by side as the first $m$ lists (i.e. the
new list in the $i$-th place is $\bigcup_{1 \leq j \leq k} l_{ji}$
), and add $m' - m$ lists disjoint from all others.

Let $T = T(L_0)$ (the transversal hypergraph of $L_0$).

Let $\{ V_\alpha \}, \alpha \in A$, be a partition of the colors in
$L_0$ into disjoint sets of colors that appear only together in the
lists of $L_0$ (i.e. if $v_1, v_2 \in V_\alpha$ then for every $l
\in L$, $v_1 \in l$ if and only if $v_2 \in l$). Then by Lemma
\ref{le4} the vertices in each $V_\alpha$ are min-equal. When
speaking of an ${R_m}$-cover, all vertices are $R_m$-equal (since
the sets in $R_m$ are all subsets of size $m-2$ of a set $S$). Let
us take a representative $v_\alpha$ from every $V_\alpha$. Each
$v_\alpha$ is given a weight of $w(v_\alpha) = |V_\alpha|$.

According to Proposition \ref{the1} there is an ${(m-2)}$-cover
$R^{*}$ of $T' = T_{cov}(\{ v_\alpha \}_{\alpha \in A})$ (the
sub-hypergraph whose edges are the minimal edges of $T$ contained in
$\{v_\alpha\}_{\alpha \in A}$) such that the cardinality of a
minimum ${(m-2)}$-cover of $T$, $s {(m-2)}^{m-2}$, equals $\sum_{r
\in R^{*}}{\{\prod_{v \in r} w(v) \} }$.


Let us use the partition ${\{ V_\alpha \}, \alpha \in A}$, to build
a partition ${\{ V_\alpha' \}, \alpha \in A \cup \{1,...,m'-m\}}$,
of the colors in $L_0'$ to sets of colors that appear only together:
We partition the colors in the first $m$ lists to sets $\{ V_\alpha'
\}, \alpha \in A$, where $V_\alpha'$ is composed of the $k$ copies
of $V_\alpha$ in the isomorphisms between $L_0$ and $L_1,...L_k$ -
then the colors in $V_\alpha'$ appear exactly where $V_\alpha$
appeared before (i.e. if $v \in V_\alpha$, $v \in l_i$, then the $k$
copies of $v$ are in $l_i'$, and if $v \notin l_i$, none of the
copies are in $l_i'$),
therefore this is a set of colors that appear only together in the
lists of $L_0'$, now of cardinality $k |V_\alpha|$. Since the
original $V_\alpha$'s were disjoint, so are the $V_\alpha'$'s. Our
last step is to add the $m' - m$ lists $l_i'$, $m < i \leq m'$, as
$m' - m$ sets in $\{ V_\alpha' \}$.


According to Lemma \ref{le4}, the vertices in each $\{ V_\alpha' \}$
are min-equal. Let us take a representative $v_\alpha'$ from every
$\{ V_\alpha' \}$. If we let one of the copies of $L_0$ be an exact
copy (for example, if $l_{1i} = l_i$ for $1 \leq i \leq m$) then we
can simply take $v_\alpha$ to represent $\{ V_\alpha' \}$ in the
first $m$ lists.

Let $H$ be the transversal hypergraph of $L_0'$. Let $H' =
H_{cov}(\{ v_\alpha' \}_{\alpha' \in A'})$ (the sub-hypergraph of
$H$ composed of all minimal edges $e \in H$ such that $e \cap
V_\alpha' \subseteq \{ v_\alpha' \}$ for every $\alpha \in A' $).

\begin{cla}\label{cl3}
$e' \in H'$ if and only if $e' = e \cup \bigcup_{a \in
\{1,...m'-m\}} {\{v_a'\}}$ for some $e \in T'$ (i.e. the edges of
$H'$ are the edges of $T'$ to which are added the representatives of
all the lists from $m+1$ on).
\end{cla}
\begin{proof}
Every edge in $H$ contains exactly one element of each of the lists
$l_i'$, $i > m$, so an edge whose intersection with $V_a'$ is a
subset of $\{ v_a' \}$, necessarily contains $v_a'$. So the set of
edges from which we need to choose minimal ones is the set of all
edges $e \in H$ such that $e \cap V_\alpha \subseteq \{ v_\alpha \}$
for every $\alpha \in A$, to which are added the representatives of
all the lists from $m+1$ on. But such an edge is minimal if and only
if the induced edge on the first $m$ lists is minimal (the edge sets
are isomorphic).
\end{proof}
\bigskip


Each representative $v_\alpha'$ is given a weight $w(v_\alpha') =
|V_\alpha'|$. If $v_\alpha'$ is contained in the first $m$ lists
(i.e. it belongs to a copy of $V_\alpha$), then $|V_\alpha'| = k
|V_\alpha|$. Otherwise $|V_\alpha'| = m'-2$ ($V_\alpha'$ is then a
list of the colors in $l_i'$ for some $i > m$).

According to Proposition \ref{the1}, in order to show that the
cardinality of a minimum ${(m'-2)}$-cover of $T(L_0')$ is at most $s
{(m'-2)}^{m'-2}$ it is enough to find an ${(m'-2)}$-cover $R'$ of
$H'$ for which $\sum_{r' \in R'}{\{\prod_{v \in r'} w(v) \} } = s
{(m'-2)}^{m'-2}$.

Let us use the $(m-2)$-cover $R^{*}$ of $T'$ to build such an
${(m'-2)}$-cover $R'$ of $H'$: $r' \in R'$ if and only if $r' = r
\cup \bigcup_{a \in \{1,...m'-m\}} {\{v_a'\}}$ for some $r \in
R^{*}$ (i.e. we add the representatives of the last $m'-m$ lists to
every $r \in R^{*}$).\\
This is a cover because if $e' \in H'$ then by Claim \ref{cl3} there
is an $e \in T'$ such that $e' = e \cup \bigcup_{a \in
\{1,...m'-m\}} {\{v_a'\}}$.
Since $e \in T'$ there is an $r \in R^{*}$ such that $r \subseteq
e$, and since $\bigcup_{a \in \{1,...m'-m\}} {\{v_a'\}} \subseteq
e'$, $r' = r \cup \bigcup_{a \in \{1,...m'-m\}} {\{v_a'\}} \subseteq
e'$.

For every $r' \in R'$, since $r' = r \cup \bigcup_{a \in
\{1,...m'-m\}} {\{v_a'\}} \subseteq e'$ for some $r \in R^{*}$,
$\prod_{v \in r'} w(v) = (m'-2)^{m'-m} {\{\prod_{v \in r} k w(v)
\}}$. Therefore $\sum_{r' \in R'}{\{\prod_{v \in r'} w(v) \} } =
(m'-2)^{m'-m} k ^ {m-2} \sum_{r \in R_{*}}{\{\prod_{v \in r} w(v) \}
} = {(m'-2)^{m'-m} {{(\frac{m'-2}{m-2})} ^ {m-2}} s {(m-2)}^{m-2}} =
s (m'-2)^{m'-2}$.

This ends the proof of Theorem \ref{thm1}.

\end{proof}

Two immediate conclusions from Theorem \ref{thm1} are:
\begin{itemize}
\item[(1)] For every even
$m$, $n_m \leq \frac{1}{2} {(m-2)}^{m-2}$. This is because $n_4 = 2$
($ch{(K_{4,2})} = 3$ and $ch{(K_{4,1})} = 2$).
\item[(2)] For every $m$ such that $3 | (m-2)$, i.e., $m \mod 3
\equiv 2$, $n_m \leq \frac{13}{27} {(m-2)}^{m-2}$. This is because
$n_5 = 13$, as shown by F\"{u}redi, Shende and Tesman in \cite{[4]}
(the configuration
$\{\{1,6,7\},\{1,8,9\},\{2,8,7\},\{2,6,9\},\{3,4,5\}\}$ with a
minimum $3$-cover of its transversal graph shows $n_5 \leq 13$).
\end{itemize}

\section{Proof of Theorem 2}


We have found a list of six $4$-tuples with a cover of size 123,
which shows $n_6 \leq 123$.

This list $L_6$ is: $\{\{1,3,5,13\},\{1,4,6,14\},\{2, 3, 7,
15\},\{2, 4, 8, 16 \},$\\
$\{5, 6, 7, 8\}, \{9,10,11,12\}\}$.

Let us generalize this structure of lists and apply it to all values
of $m$, to show:

\begin{thm}\label{thm2} There is an infinite sequence of integers $S$, such that if ${m
\in S}$, then $n_m \leq 0.4643 {(m-2)}^{m-2}$.
\end{thm}

\begin{proof}
Let $k_1, k_2, k_3, k_4$, and $k_5$ be integers such that $\sum_{1
\leq i \leq 5}{k_i} = m-2$, $4k_4 \leq m-2$, $0 \leq k_i \leq m-2$
for every $1 \leq i \leq 5$. Let $L$ be a set of $m$ $(m-2)$-tuples,
divided into five lists which will be described next, and $m-5$
lists that are disjoint from all other lists. The first five lists
$l_1,...,l_5$ are built as follows: the intersection of every
threesome of these lists is empty, $|l_1 \cap l_2| = |l_3 \cap l_4|
= k_1$, $|l_1 \cap l_3| = |l_2 \cap l_4| = k_2$, $|l_1 \cap l_4| =
|l_2 \cap l_3| = k_3$, and $|l_1 \cap l_5| = |l_2 \cap l_5| = |l_3
\cap l_5| = |l_4 \cap l_5| = k_4$.


We can partition the colors of $L$ into disjoint sets according to
the lists in which they appear together as follows: $V_1$ = $\{l_1
\cap l_2\}$, $V_2$ = $\{l_3 \cap l_4\}$, ${V_3 = \{l_1 \cap l_3\}}$,
$V_4$ = $\{l_2 \cap l_4\}$, $V_5$ = $\{l_1 \cap l_4\}$, $V_6$ =
$\{l_2 \cap l_3\}$, $V_7$ = $\{l_1 \cap l_5\}$, ${V_8 = \{l_2 \cap
l_5\}}$, $V_9$ = $\{l_3 \cap l_5\}$, $V_{10}$ = $\{l_4 \cap l_5\}$,\\
$V_{10+i}$ = $l_i \setminus \bigcup_{1 \leq j \leq 5, j \neq
i}{\{l_i \cap l_j\}}$ ($|V_{10+i}| = k_5$) for $1 \leq i \leq 4$,\\
$V_{15}$ = $l_5 \setminus  \bigcup_{1 \leq j \leq 4, }{\{l_5 \cap
l_j\}}$ ($|V_{15}| = m-2-4{k_4}$), and\\ $V_{15+i}$ = $l_{5+i}$
($|V_{15+i}| = m-2$) for $1 \leq i \leq m-5$.

According to Proposition \ref{the1}, the cardinality of the smallest
$(m-2)$-cover of $T(L)$ is equal to the the value of a minimum
weighted $(m-2)$-cover of the transversal hypergraph of a set of
lists in which we take one representative $v_i$ from each of the
$V_i$'s, and give $v_i$ weight $|V_i|$.

Taking one representative from each $V_i$, gives us a set of $m$
lists $L'$ with the following structure:\\
1  3  5  7  11 \\
1  4  6  8  12 \\
2  3  6  9  13 \\
2  4  5  10  14 \\
7  8  9  10  15 \\
16 \\
17 \\
...\\
with weights $w(1) = w(2) = k_1$, $w(3) = w(4) = k_2$, $w(5) = w(6)
= k_3$, $w(7) = w(8) = w(9) = w(10) = k_4$, $w(11) = w(12) = w(13) =
w(14) = k_5$, $w(15) = m-2-4k_{5}$, and $w(16) = w(17) = ... = m-2$
(there are $m-5$ vertices of this last type).

If we set ${\alpha}_i = \frac{k_i}{m-2}$, then the weight of any
$(m-2)$-tuple, divided by ${(m-2)}^{m-2}$, is a function of the
${\alpha}_i$'s alone - let us from now on omit ${(m-2)}^{m-2}$.

Let us describe a cover for the above hypergraph (it is possible to
prove that this is a minimal cover of this structure, but we shall
not prove it in this paper): Start with the set $T$ of all tracks of
$L'$. At step $i$ take all the transversals of cardinality $m-2$
that belong to at least one track in $T$ as edges in the cover.
Remove the tracks that belong to these transversals from the set
$T$. Now all the tracks in $T$ belong to transversals of cardinality
$m-1$ or more, so we remove the last coordinate from every track in
$T$, and continue to step $i+1$. We end the process when $T$ is
empty (after two steps all the tracks left have $m-2$ vertices).

A straightforward calculation shows that the value (sum of weights)
of this $(m-2)$-cover is: ${\alpha_1}^2 + {\alpha_2}^2 +
{\alpha_3}^2 + 2({{{\alpha}_4}^2 + 2{{\alpha}_4}{{\alpha}_5}} + {(1
- 2{{\alpha}_4}){{{\alpha}_5}^2}})({\alpha}_1 + {\alpha}_2 +
{\alpha}_3) +\\ {4{({{\alpha}_4} + {{\alpha} _5}{(1 -
3{{\alpha}_4})})}{({{\alpha}_1}{{\alpha}_2} +
{{\alpha}_1}{{\alpha}_3} + {{\alpha}_2}{{\alpha}_3})}} +
4{{\alpha}_1}{{\alpha}_2}{{\alpha}_3}{(1 - 3{{\alpha}_4})} +
({({{\alpha}_4} + {{\alpha}_5})}^4 - {{\alpha}_5}^4 ) +
{{\alpha}_5}^4 ( {1 - 4{{\alpha}_4}} )$.

Putting in ${{\alpha}_1} = {{\alpha}_2} = {{\alpha}_3} =
\frac{1}{3}$, ${{\alpha}_4} = {{\alpha}_5} = 0$ gives
$\frac{13}{27}$, which we know is the value for $m = 5$.

Putting in ${{\alpha}_1} = {{\alpha}_2} = {{\alpha}_4} =
{{\alpha}_5} = \frac{1}{4}$, ${{\alpha}_3} = 0$ gives
$\frac{123}{256}$, which is the value we calculated for $m = 6$
(this is better than Eaton's result of $\frac{125}{256}$, mentioned
in Tuza's survey paper \cite{[5]}).

Now we wish to see when this is brought to a minimum as a function
of ${{\alpha}_1}, {{\alpha}_2}, {{\alpha}_3}, {{\alpha}_4},
{{\alpha}_5}$.  For a given $m$ we can only take ${{\alpha}_i}$'s
such that ${{\alpha}_i} {(m-2)}$ is an integer for $1 \leq i \leq 5$
- the result we get will be applicable only to such $m$'s.

Minimizing by setting ${{\alpha}_3} = {{\alpha}_1} + \varepsilon$
and equating the derivative of the expression as a function of
$\varepsilon$ to $0$, gives ${{\alpha}_1} = {{\alpha}_2}
={{\alpha}_3}$ or ${{\alpha}_1} = {{\alpha}_2}$, ${{\alpha}_3} = 0$.

In the case ${{\alpha}_1} = {{\alpha}_2} ={{\alpha}_3}$ we get that
the minimum cover is of cardinality $\frac{13}{27} - \frac{6}{27}
{{\alpha}_4} + \frac{13}{9} {{\alpha}_4}^2 - \frac{58}{27}
{{\alpha}_4}^3 + \frac{13}{9} {{\alpha}_4}^4 + {{\alpha}_5} (
\frac{2}{9} -\frac{22}{9} {{\alpha}_4} + \frac{26}{9} {{\alpha}_4}^2
+ \frac{4}{3} {{\alpha}_4}^3 ) + {{{{\alpha}_5}^2} ( \frac{1}{9}
-\frac{2}{9} {{\alpha}_4} + \frac{10}{3} {{\alpha}_4}^2 )} +
{{{\alpha}_5}^3} ( -\frac{22}{27} + \frac{40}{9} {{\alpha}_4} ) +
{{{\alpha}_5}^4} (1 - 4{{\alpha}_4} )$.

Numerical optimization gives a value of 0.4642..., for ${\alpha}_4 =
0.1969..., {\alpha}_5 = 0.2123...$.

\end{proof}

\newpage

\section{Proof of Theorem 3}


\begin{thm}\label{thm3}
If $m \rightarrow \infty$, then $n_m$ is asymptotically at most
$0.474 {(m-2)}^{m-2}$.
\end{thm}

\begin{proof}
Let us take the structure of $L$ in which $|l_1 \cap l_2 \cap l_3| =
k \geq 2$, ($|l_1 \cap l_2 \cap l_3 \cap l_i| = 0$ for $4 \leq i
\leq m-k+1$), and $|l_1 \cap l_{i}| = 1$ for $4 \leq i \leq m-k+2$.
Also $|(l_2 \cap l_3) \setminus l_1| = l$. All the rest of the
colors are
different:\\
1 2 3... 4 5 ... 6  7...\\
1 2 3... 8 9 ... 10 11...\\
1 2 3... 8 9 ... 12 13...\\
4 14...\\
5 15...\\
...\\
6 17...\\
7 18...\\
...\\
19...\\
...\\

The cover in this case is: All $(m-2)$-tuples that represent tracks
in which the same color $c$ in $l_1 \cap l_2 \cap l_3$ is chosen
from $l_1, l_2, l_3$ - there are $k {(m-2)}^{m-3}$ such
$(m-2)$-tuples; All $(m-2)$-tuples that represent tracks in which
the same color $j \in (l_2 \cap l_3) \setminus l_1$ is chosen out of
lists $l_2$ and $l_3$, and the color chosen out of $l_1$ is chosen
again out of the other list in $l_4,...$ which contains it - there
are $l ({(m-2)^{m-2-k} - {(m-3)}^{m-2-k}}){(m-2)^{k-1}}$ such
tracks; All the tracks that have not been covered yet belong to
transversals of cardinality $m-1$ - so we remove the last coordinate
in every track, and take all the minimal transversals that are left
- each belongs to exactly one track in which the color chosen out of
$l_1$ is also chosen by the other list that contains it - there are
${(m-2-k-l)}^2 ({(m-2)}^{m-2-k}-{(m-3)}^{m-2-k}) {(m-2)}^{k-2}$ such
tracks.

All in all the cover has cardinality\\ $k {(m-2)}^{m-3} + l
({(m-2)^{m-2-k} - {(m-3)}^{m-2-k}}){(m-2)^{k-1}} + {(m-2-k-l)}^2
({(m-2)}^{m-2-k}-{(m-3)}^{m-2-k}) {(m-2)}^{k-2} =  k {(m-2)}^{m-3} +
{(m-2)}^{k-2}[{ {(m-2)}^{m-2-k}-{(m-3)}^{m-2-k} ]} (
{(m-k-2)}^2 + l^2 +2lk - l (m-2) )$.\\
and this expression is minimal as a function of $l$ when $l =
\frac{m-2-2k}{2}$ (since $l$ must be an integer, if $m$ is odd we
will take $l = \frac{m-3-2k}{2}$, and the rest of the proof is
similar). Putting this value of $l$ back into the
expression gives \\
$k {(m-2)}^{m-3} + \\ {(m-2)}^{k-2}{[
{(m-2)}^{m-2-k}-{(m-3)}^{m-2-k} ]} ( \frac{3}{4} {{(m-2)}^2} - k
{(m-2)} )$.

When $m \rightarrow \infty$, the expression tends to \\ $ k
{(m-2)}^{m-3} + {(m-2)}^{k-2}{[ {(m-2)}^{m-2-k} {(1-\frac{1}{e})} ]}
( \frac{3}{4} {{(m-2)}^2} - k {(m-2)} ) =
 ( { \frac{k}{e (m-2)} + \frac{3}{4}{(1-\frac{1}{e})}}
) {(m-2)}^{m-2} $.

This expression is minimal when $k$ is as small as possible - in
this case when $k = 2$. If this is the case we get that the
cardinality of the cover tends to $ ( {
\frac{3}{4}{(1-\frac{1}{e})}} ){(m-2)}^{m-2} = 0.474 {(m-2)}^{m-2}$.

\end{proof}
\newpage

\section{Conclusion and Open Problems}
In this paper we used certain structures of the family of lists
assigned to $M$ ($|M| = m$) to calculate upper bounds on $n_m$, the
smallest integer $n$ such that $ch(K_{m,n}) = m-1$. These are not,
of course, all the possible structures of a family of lists on $M$
with all transversals of cardinality at least $m-2$ (the transversal
hypergraph of a family of $m$ lists which contain a transversal of
cardinality less than $m-2$ does not allow an $(m-2)$-cover and
therefore does not need to be considered), and the remaining
structures still need to be analyzed. The method we have devised for
finding a minimum cover of a hypergraph by another hypergraph, of
solving an equivalent problem of finding a minimum weighted cover,
can also be used to calculate upper bounds on the smallest integer
$n$ such that $ch(K_{m,n}) = m-k$ for other small $k$'s.

\newpage

\end{document}